\documentclass[letterpaper, 10 pt, conference]{ieeeconf}  % Comment this line out
\IEEEoverridecommandlockouts                              
\usepackage{cite}

\usepackage{amsthm}
\usepackage{amsmath,amssymb,amsfonts}
\usepackage{algorithmic}
\usepackage{graphicx}
\usepackage{algorithm,algorithmic}
\usepackage{hyperref}
\usepackage{amsmath,amssymb,amsfonts}
\usepackage{bm}
\usepackage{algorithm}
\usepackage{graphics}
\usepackage{amsfonts}
\usepackage{amssymb}
\usepackage{mathrsfs}%����

\usepackage{graphics} % for pdf, bitmapped graphics files
\usepackage{epsfig} % for postscript graphics files
\usepackage{amsmath} % assumes amsmath package installed
\usepackage{amssymb}  % assumes amsmath package installed
\usepackage[dvipsnames]{xcolor}

\newtheorem{theorem}{Theorem}[section]
\newtheorem{remark}[theorem]{Remark}
\newtheorem{definition}[theorem]{Definition}
\newtheorem{lemma}[theorem]{Lemma}

\newtheorem{assumption}[theorem]{Assumption}

\title{ Optimal Adaptive Control of  Linear Stochastic Systems with Quadratic Cost Function}
\author{Nian Liu, Cheng Zhao, Shaolin Tan, and Jinhu L{\"{u}}, \emph{Fellow, IEEE}
}

\begin{document}

\maketitle
\thispagestyle{empty}
\pagestyle{empty}

\begin{abstract}
	In this paper, we consider the adaptive linear quadratic Gaussian control problem, where both the linear transformation matrix of the state $A$ and the control gain matrix $B$ are unknown. The proposed adaptive optimal control only assumes that $(A, B)$ is stabilizable and $(A, Q^{1/2})$ is detectable, where $Q$ is the weighting matrix of the state in the quadratic cost function. This condition significantly weakens the classic assumptions used in the literature. To tackle this problem, a weighted least squares algorithm is modified by using  random regularization method, which can ensure uniform stabilizability and uniform detectability of the family of estimated models. At the same time,  a diminishing excitation is incorporated into the design of the proposed adaptive control to guarantee strong consistency of the desired components of the estimates.  Finally, by utilizing this family of estimates, even if not all components of them converge to the true values,  it is demonstrated that a certainty equivalence control with such a diminishing excitation  is optimal for an ergodic quadratic cost function. 
\end{abstract}

\section{Introduction}

\label{sec:introduction}
Learning from data is pivotal in advancing scientific knowledge and understanding across diverse disciplines. It serves as the bedrock of artificial intelligence and is increasingly prevalent, particularly in control engineering. Adaptive control stands out as a classical theory rooted in learning from data, with extensive research dedicated to it from the end of 1950s (see, e.g., \cite{a3,b1,a4}). Over the past six decades, substantial progress has been achieved in this field (see, e.g., \cite{b2,b3,a5}).

Adaptive control is a specialized area within control engineering that focuses on designing control systems capable of adjusting their parameters and behavior in response to changes in the system or its environment.  The primary goal of adaptive control is to achieve stable and optimal performance in the presence of uncertainties, variations, and disturbances, which may be difficult to model or predict accurately. One prominent approach in adaptive control design is the certainty equivalent principle, which involves using observed information to estimate unknown parameters at each time instant and then updating the controller by taking the estimates as the “true” parameter. While this method may not be optimal at every time instant, it holds strong appeal in practical design and proves effective in addressing dynamical systems with unknown parameters. Nevertheless, establishing a rigorous theory on the stability and convergence of stochastic adaptive control, even for linear stochastic uncertain systems, presents considerable challenges. This difficulty arises from the complex nonlinear stochastic dynamics that typically describe closed-loop control systems in such scenarios.

For example, one of the most well-known longstanding problem in adaptive control was how to establish a rigorous convergence theory for the classical self-tuning regulators, which was solved by a natural combination of the least-squares estimate with the minimum variance control for linear stochastic systems (see, e.g., \cite{a6,a7}). Another  well-known example was the optimal adaptive linear quadratic control problem, where a key theoretical difficulty was how to guarantee both the global stability of the closed-loop adaptive systems and the controllability of estimated models, which was initially solved for general open-loop stable linear stochastic systems in \cite{a8,a10}.
Without resorting to such an open-loop stability condition, this longstanding problem was solved in the work \cite{a1} for discrete-time linear systems, based on the self-convergence property of weighted least squares (WLS) algorithm. As for adaptive continuous-time linear quadratic Gaussian (LQG) control problem, under the assumption that $(A, B)$ is controllable and $(A,Q^{1/2})$ is observable, it was later reasonably solved in the work \cite{a2}  based on the self-convergence property of continuous-time WLS algorithm and on a random regularization theory established in \cite{a1}. The above works turn out to be directly instrumental in solving the adaptive problems proposed in this paper under the conditions that $(A, B)$ is stabilizable and $(A, Q^{1/2})$ is detectable.

%% data-driven and reinforcement learning
In a related but somewhat different domain, data-driven control aims to design controllers directly from data, bypassing the explicit identification of a system model. This approach is not only conceptually appealing but also valuable in scenarios where system identification is challenging or even infeasible due to insufficient data. A seminal contribution to data-driven control is often credited to Ziegler and Nichols \cite{d1} for their initial work on tuning proportional–integral–derivative controllers. More recently, significant attention has been devoted to the problem of deriving optimal controllers from data \cite{d2,d3}. Beyond control problems, \cite{d4} analyzes the stability of an input/output system using time series data and the papers \cite{d5,d6,d7} deal with data-based controllability and observability analysis. 

Besides, reinforcement learning (RL) techniques have garnered considerable attention for addressing optimal control problems across both discrete and continuous-time frameworks. For example, \cite{r1} and \cite{r2} study the linear quadratic (LQ) regulator control problems in completely deterministic scenarios. Global convergence of policy gradient methods for deterministic LQ regulator problem is shown in \cite{r3}. In \cite{r4}, the authors tackle the stochastic LQ problem under the assumption that the noise is measurable, and propose a least squares temporal difference learning approach.
%In \cite{r5}, the authors give LSTD algorithms for the LQ problem with process noise and analyze the regret bound, i.e., the difference between the occurred cost and the optimal cost. 
%The regret bound is also studied in \cite{r6} for a model-building routine. The sample complexity of the LSTD for the LQ problem is studied in \cite{r7}. \cite{r8} study the effect of measurement noise in the proposed algorithms and the classical off-policy and the classical Q-learning routines. Bradtke et al. 
\cite{r9} introduced a Q-learning policy iteration method aimed for addressing discrete-time LQ problems. This technique relies on the concept of the Q-function, as elucidated in \cite{r11}.
%For its recent applications to discrete-time models, we refer to Rizvi and Lin \cite{r12}, Luo et al. \cite{r13}, and Kiumarsia et al. \cite{r14}. Baird \cite{r15} first adopted an RL approach to obtain the optimal control for a continuous-time discrete-state system. 
%Murray et al. \cite{r16} proposed an iterative adaptive dynamic programming (ADP) scheme for nonlinear systems. 
Recently, a number of new RL methods
have been developed for optimal control problems in continuous-time
cases (e.g., \cite{r16,r17,r18,r19}).  

In this paper, a complete solution to the continuous-time adaptive LQG control problem is given, using only the natural
assumptions of stabilizability and detectability for the system matrices. The optimal adaptive control is designed as following steps. The initial step involves employing the WLS algorithm to generate a family of convergent estimates, without imposing constraints on excitation or closed-loop system stability, demonstrating the self-convergence property introduced in \cite{a2}. Subsequently, the estimates are refined through a random regularization procedure outlined in \cite{a1}, resulting in a uniformly stabilizable and detectable family of estimated models. Additionally, integrating a diminishing excitation (white noise) into the certainty equivalence control ensures the strong consistency of some desired components of the estimates, and at the same time, such adaptive control optimizes the quadratic cost function. The main contributions of this paper are summarized as follows: 
\begin{enumerate}
	\item  This paper gives a concrete approach for the design of the optimal adaptive control under the assumption that $(A, B)$ is stabilizable and $(A, Q^{1/2})$ is detectable, which greatly weakens the conditions previously used in \cite{a1} and \cite{a2}.
	\item This paper initially introduces an original concept of uniformly stablizability and uniformly detectability, which can be guaranteed for the estimated models by employing the random regularization method and is useful for the design of adaptive control.
	\item  This paper shows that if the excitation associated with some components of parameter estimates is rich, accurate estimation of these parameters can be obtained, and they are unaffected by the estimates derived from insufficient excitation. We emphasize that even if the estimates may not converge to the true value, the proposed adaptive control still has the  ability to achieve global stability and optimality.
\end{enumerate}

The remainder of the paper is organized as follows: In Section \ref{section2}, we introduce the adaptive LQG problem, and the WLS estimator associated with the random regularizatin method. In Section \ref{section3}, the design procedure of the optimal adaptive control is presented and the main results on stability and
optimality of the closed-loop systems are established. Section \ref{section4} gives the proof of the main theorems. Finally, some concluding remarks are given in Section \ref{section6} and some proofs of the lemmas are given in Appendix.

\section{Problem Formulation}\label{section2}
\subsection{Notations}
In this paper, we use $ A\in\mathbb{R}^{m\times n}$ to denote an  $m\times n$-dimensional real matrix. $ I_n$ is the identity matrix of $\mathbb{R}^{n\times n}$. For a matrix $ A\in\mathbb{R}^{n\times n}$, we use $\lambda_{\min}(A)$ to denote the minimum eigenvalue of $A$, while $\lambda_{\max}(A)$ denotes the maximum eigenvalue of $A$. $\| A\|$ denotes the Euclidean norm of $A$, i.e., $\| A\|=(\lambda_{\max}( A A^{\mathsf{T}}))^{\frac{1}{2}}$, where ${\mathsf{T}}$ denotes the transpose operator. In addition, $|\cdot|$ denotes the vector $L_2$-norm and $\langle  \cdot,  \cdot\rangle$ denotes the inner product of vectors. $A >0$ ($ A\ge0$) means that $A$ is a positive definite (semi-definite) matrix. $\text{tr}(A)$ denotes the trace of $A$. Let $\lfloor t\rfloor$ represent the integer part of $ t$ and $\mathbb{N}$ be the set of natural numbers. Let $\sigma(A)$ be the set of all eigenvalues of the matrix $A$, $\sigma^+(A)\triangleq\{\lambda\in\sigma(A)~\!|~\!\text{Re}(\lambda)\ge0\}$, and $\sigma^-(A)\triangleq\{\lambda\in\sigma(A)~\!|~\!\text{Re}(\lambda)<0\}$. For two matrices $A\in\mathbb{R}^{n\times n}$, $B\in\mathbb{R}^{n\times m}$, let us denote $\text{Im}(B)=\{B\eta:\eta\in\mathbb{R}^m\},$  $\text{Im}(A|B)= \text{Im}(B)+\cdots+A^{n-1}\text{Im}(B),$ where $A~\!\text{Im}(B)=\{A\eta:\eta\in \text{Im}(B)\}$.

\subsection{LQG Control}

Consider the following LQG control systems:
\begin{equation}\label{s1}
	\mathsf{d} x(t)=(Ax(t)+Bu(t))\mathsf{d} t+D\mathsf{d} w(t),
\end{equation}
where $x(t)\in \mathbb{R}^n$ is the state of the system with initial condition $x(0)$, $u(t)\in\mathbb{R}^{m}$ is the control, and $(w(t),\mathscr{F}_t;t\ge0)$ is a $\mathbb{R}^p$-valued standard Wiener process. $A\in\mathbb{R}^{n\times n}, B\in\mathbb{R}^{n\times m},D\in \mathbb{R}^{n\times p}$ are \emph{unknown} matrices.

The objective is to design an admissible control to minimize the following ergodic cost function:
\begin{equation}\label{cost}
	J(u)=\limsup\limits_{T\to \infty}\frac{1}{T}\int_{0}^{T}(x^{\mathsf{T}}Qx+u^{\mathsf{T}}Ru)\mathsf{d} t,
\end{equation}
where $Q\ge0$ and $R>0$ are known weighting matrices.  

First, we introduce the definition of admissible control.
\begin{definition}\cite{a2}
	A control $\{u(t),t\ge0\}$ for the system (\ref{s1}) is said to be admissible if it is adapted to $\{\mathscr{F}_t;t\ge0\}$ and under which, the following properties hold for any initial state $x(0)\in\mathbb{R}^n$:
	\begin{equation}\notag
		\limsup\limits_{T\to \infty}\frac{\int_{0}^{T}|u(t)|^2\mathsf{d}t}{\int_{0}^{T}|x(t)|^2\mathsf{d}t}<\infty\quad a.s.
	\end{equation}
\end{definition}

\begin{remark}
	Here, $\{\mathscr{F}_t;t\ge0\}$ is a filtration containing the filtration  $\mathscr{F}^x_t\triangleq \sigma(x(s);0\le s\le t)$ generated by the state. Such definition can make the design of adaptive control to be more flexible to include, e.g., diminishing excitation or “exploration” signals.
\end{remark}

It is well-known that if the matrices $A,B$ are known, then  the following algebraic Riccati equation (ARE):
\begin{equation}\label{are1}
	A^{\mathsf{T}}X+XA+Q-XBR^{-1}B^{\mathsf{T}}X=0
\end{equation}
admits a stabilizing solution $X$, i.e., $X\ge0$ and $A-BR^{-1}B^{\mathsf{T}}X$ is stable. Moreover, the feedback control 
\begin{equation}\label{u}
	u(t)=-R^{-1}B^{\mathsf{T}}Xx(t)
\end{equation} can stabilize the system and minimize the cost function with 
\begin{equation}
	J(u)=\text{tr}(D^{\mathsf{T}}XD) \quad a.s.
\end{equation}
We remark that if $(A,B)$ is stabilizable and $(A,Q^{1/2})$ is detectable,  \textcolor{blue}{if and only if}  such stabilizing solution exists and it is unique.

However, if the system matrices $A,B$ are unknown, then the optimal feedback control (\ref{u}) is not available. Therefore, it is natural to consider the adaptive control. It is worth mentioning that, the papers \cite{a1} and \cite{a2} designed the adaptive control  for continuous-time and discrete-time systems respectively, under the conditions that $(A,B)$ is controllable and $(A,Q^{1/2})$ is observable. But, such conditions on matrices $\{A,B,Q\}$ are somewhat restrictive and may not be satisfied by some practical situations. A natural question is, can we construct an optimal adaptive control under the following weaker assumption.

\begin{assumption}\label{ass1} 
	The pair $(A,B)$ is stabilizable and the pair $(A,Q^{1/2})$ is detectable.
\end{assumption}

In this paper, we will prove that Assumption \ref{ass1} is \emph{sufficient} for the existence of an optimal adaptive control.
To give a concrete construction, we first introduce WLS estimator and random regularization method.

\subsection{The WLS Estimation} 
To describe the adaptive problem in the standard form, we  introduce the following notations:
\begin{equation}\label{key9}
	\theta ^{\mathsf{T}}=[A,B]
\end{equation}
and
\begin{equation}
	\varphi(t)=[x^{\mathsf{T}}(t),u^{\mathsf{T}}(t)]^{\mathsf{T}},
\end{equation}
and rewrite the system (\ref{s1}) into the following compact form:
\begin{equation}
	\mathsf{d}x(t)=\theta^{\mathsf{T}}\varphi(t)\mathsf{d}t+D\mathsf{d}w(t).
\end{equation}
The continuous-time WLS estimates $\{\theta(t),t\ge0\}$ are given by \cite{a2}
\begin{equation}\label{wls1}
	\mathsf{d}\theta(t)=a(t)P(t)\varphi(t)[\mathsf{d}x^{\mathsf{T}}(t)-\varphi^{\mathsf{T}}(t)\theta(t)\mathsf{d}t],
\end{equation}
and
\begin{equation}\label{wls2}
	\mathsf{d}P(t)=-a(t)P(t)\varphi(t)\varphi^{\mathsf{T}}(t)P(t)\mathsf{d}t,
\end{equation}
where $P(0)=I_{n+m}$, $\theta^{\mathsf{T}}(0)=[A(0),B(0)]$ are arbitrary deterministic values such that $(A(0),B(0))$ is stabilizable and $(A(0),Q^{1/2})$ is detectable, and
\begin{equation}\label{a1}
	a(t)=\frac{1}{\log^2r(t)},
\end{equation}
with
\begin{equation}\label{r}
	r(t)=\|P^{-1}(0)\|+\int_{0}^{t}|\varphi(s)|^2\mathsf{d}s,
\end{equation}
It can be verified that $P(t)$ defined in (\ref{wls2}) has the following explicit form:
\begin{equation}\label{p}
	P(t)=\Big(P^{-1}(0)+\int_{0}^{t}a(s)\varphi(s)\varphi^{\mathsf{T}}(s)\mathsf{d}s\Big)^{-1}.
\end{equation}
The nice properties of WLS estimates are presented in the following lemma.
\begin{lemma} \cite{a2}\label{le1}
	Let $\{\theta(t),t \ge 0\}$ be defined by (\ref{wls1})-(\ref{r}). Then the
	following properties are satisfied:
	\begin{flalign}
		&1) \sup\limits_{t\ge0}\|P^{-\frac{1}{2}}(t)(\theta(t)-\theta)\|^2<\infty\quad a.s.\notag&&\\
		&2)\int_{0}^{\infty}a(t)\|(\theta(t)-\theta)^{\mathsf{T}}\varphi(t)\|^2 \mathsf{d} t<\infty\quad a.s.\notag&&\\
		&3)\lim\limits_{t\to\infty}\theta(t)=\bar{\theta}\quad a.s.\notag&&
	\end{flalign}
	where $\bar{\theta}$ is a certain random matrix.
\end{lemma}
\begin{remark}
	From Lemma \ref{le1}, one can see that the WLS estimates are self-convergent, but it should be noted that the family of estimated models (\ref{wls1})-(\ref{r}) may not be stabilizable and detectbable even if Assumption \ref{ass1} is satisfied. This is one of the main difficulties in
	the adaptive LQG control problem considered in this paper, which could be overcome by using the method of random regularization method introduced in \cite{a1}.
\end{remark}
Before introducing the method of random regularization, let us denote the following two matrix functions:
$$Z(s,A,B)\!=\!\left[\left(\begin{array}{cc}
	s_1I_n-A & s_2I_n  \\ 
	-s_2I_n & s_1I_n-A \\ 
\end{array}\right),\left(\begin{array}{cc}
	B  & 0 \\ 
	0 & B  \\ 
\end{array}\right)\right],$$
where $s=s_1+s_2j$ and $j$ is the imaginary unit, and
\begin{align}\label{y}
	Y(A,B)&=\prod_{s\in\sigma^+(A)}\det\left(Z(s,A,B)Z^{\mathsf{T}}(s,A,B)\right) \notag\\
	&\quad\cdot\prod_{s\in\sigma^-(A)}\det\left(Z(-s,A,B)Z^{\mathsf{T}}(-s,A,B)\right),
\end{align}
where $\sigma^+(A)\triangleq\{\lambda\in\sigma(A)~\!|~\!\text{Re}(\lambda)\ge0\}$, and $\sigma^-(A)\triangleq\{\lambda\in\sigma(A)~\!|~\!\text{Re}(\lambda)<0\}$.
With the above notations, we have the following lemma and the proof is given in Appendix.
\begin{lemma}\label{le2}
	The following statements are equivalent:\\
	1) $(A,B)$ is stabilizable;\\
	2) $(sI-A,B)$ is full of row rank for any $s\in\sigma^+(A)$;\\
	3) $Z(s,A,B)$ is full of row rank for any $s\in\sigma^+(A)$;\\
	4) $Y(A,B)>0$.
\end{lemma}

From Lemma \ref{le2}, we introduce the following concept of uniform stabilizability (detectability), which is an extension of the uniform controllability (observability) defined in \cite{a1}, and will be useful in the  design and analysis of the proposed adaptive control.
\begin{definition}[Uniformly stabilizable (detectable)]
	A family of matrix pairs $\{A(t),B(t); t \ge 0\}$ with $A(t)\in \mathbb{R}^{n\times n},B(t) \in \mathbb{R}^{n\times m}$ is said to be uniformly stabilizable if there is a constant $c > 0$ such that
	$$Y(A(t),B(t))\ge c,\text{ for all } t \in [0,\infty).$$
	A family of matrix pairs $\{A(t),C(t); t \ge 0\}$ with $C(t)\in \mathbb{R}^{p\times n}$ is said to be uniformly detectable if $\{A^T(t),C^T(t); t \ge 0\}$ is uniformly stabilizable.
\end{definition}

\subsection{Random Regularization}
By Lemma \ref{le1} 3), we know that $\theta(t)$ converges to a certain random matrix $\bar{\theta}$, which may not be the true parameter matrix $\theta$. Therefore, the stabilizability and detectability properties of the corresponding estimated models $(A(t),B(t))$ and  $(A(t),Q^{1/2})$ needed for constructing the adaptive version of the control may not be guaranteed. In view of this, we resort to the regularization method introduced in \cite{a1} to modify the WLS estimates, so that the family of estimates $(A(t),B(t))$ is uniformly uniform stabilizable and  $(A(t),Q^{1/2})$ is uniformly detectable. 

To address this issue, from the boundedness of the matrix sequence $\{P^{-\frac{1}{2}}(t)(\theta-\theta(t)),\,\,t\ge0\}$, as established in Lemma \ref{le1} 1), there exists a bounded sequence $\{\beta^*(t),t\ge0\}$ such that
\begin{equation}\label{b}
	\theta=\theta(t)-P^{\frac{1}{2}}(t)\beta^*(t),
\end{equation}
which serves as a basis for a modification to the WLS estimates:
\begin{equation}\notag
	\theta(t,\beta(t))=\theta(t)-P^{\frac{1}{2}}(t)\beta(t),
\end{equation}
where $\{\beta(t)\in\mathbb{R}^{(m+n)\times n}, t\ge0\} $ is a sequence of bounded matrices to be chosen shortly. For simplicity, we rewrite $\theta(t,\beta(t))$ in the following form:
\begin{equation}\notag
	\theta^{\mathsf{T}}(t,\beta(t))=[\bar{A}(t),\bar{B}(t)].
\end{equation}

Then, we aim to choose suitable $\{\beta(t),t\ge0\}$ to ensure that $\{\bar{A}(t),\bar{B}(t);t\ge0\}$ is uniformly stabilizable and $\{\bar{A}(t),Q^{1/2};t\ge0\}$ is uniformly detectable.
From Lemma \ref{le2}, it suffices to select  $\{\beta(t),t\ge0\}$ to guarantee the uniform positivity of $f_t(\beta(t))\triangleq Y(\bar{A}(t),\bar{B}(t))Y(\bar{A}^{\mathsf{T}}(t),Q^{1/2})$, i.e., $f_t(\beta(t)\ge c$ for all $t\ge0$ with some constant $c>0$.
\begin{remark}
	In \cite{a1} and \cite{a2}, the authors select a proper sequence $\{\beta(t),t\ge0\}$ to guarantee that
	\begin{align}\label{f1}
		\bar{f}_t(\beta(t))&=\det\big(\sum_{i=0}^{n-1}\bar{A}^i(t)\bar{B}(t)\bar{B}^{\mathsf{T}}(t)\bar{A}^{i{\mathsf{T}}}(t)\big)\notag\\
		&\quad\cdot\det\big(\sum_{i=0}^{n-1}\bar{A}^i(t)Q^{1/2}(Q^{1/2})^{\mathsf{T}}\bar{A}^{i{\mathsf{T}}}(t)\big)\
	\end{align}
	is uniformly positive,
	which in turn ensures the uniform controllability of $(\bar{A}(t),\bar{B}(t))$ and the uniform observability of $(\bar{A}(t),Q^{1/2})$.
	We would like to mention that
	the construction of the above   function $f_t(\beta(t))$ is inspired by $\bar{f}_t(\beta(t))$ given by (\ref{f1}). In fact, it is not difficult to see that, if $Y(A,B)\triangleq\prod_{s\in\sigma(A)}\det\left(Z(s,A,B)Z^{\mathsf{T}}(s,A,B)\right)$, the positivity of $f_t(\beta(t))$ is equivalent to the positivity of $\bar{f}_t(\beta(t))$.
	%In fact, the construction of $Y(A,B)$ follows the following process. First, we recall that $(A,B)$ is stabilizable if and only if $(sI-A,B)$ is full of row rank for any $s\in\sigma^+(A)$, which makes us define $Y_0(A,B)$ as
	%$$Y_0(A,B)=\prod_{s\in\sigma^+(A)}\det\left([sI-A,B][sI-A,B]^{\mathsf{T}}\right).$$
	%In order to guarantee $Y_0(A,B)\in \mathbb{R}$, we have noticed Lemma \ref{le2}, so $Y_0(A,B)$ may be defined by 
	%$$Y(A,B)=\prod_{s\in\sigma^+(A)}\det\left(Z(s,A,B)Z(s,A,B)^{\mathsf{T}}\right).$$
	%But, it is not continuous respect to $A$. Finally, to ensure the continuity of $Y_0(A,B)$, we have found the desired form as (\ref{y}).
\end{remark}	

Next, we give a concrete construction of the sequence $\{\beta(t),t\ge0\}$. Let $\{\eta_k\in\mathbb{R}^{(n+m)\times n}, k\in\mathbb{N}\}$ be a sequence of random variables that are uniformly distributed in the unit ball of $\mathbb{R}^{(m+n)\times n}$ and are also independent of $\{w(t),t\ge0\}$. The sequence $\beta(k)$ is recursively given by the following:
\begin{align}
	&\beta(0)=0,\notag \\
	&\beta(k)=\begin{cases}
		\eta_k,&\text{if } f_k(\eta_k))\ge\gamma f_k(\beta(k\!-\!1)),\\
		\beta(k \!-\!1),&\text{otherwise},
	\end{cases}\notag
\end{align}
where $\gamma\in(1,\sqrt{2})$ is a fixed constant. Thus, a sequence of
regularized estimates $\{\bar{\theta}(k), k\in\mathbb{N}\}$ can be defined by 
\begin{equation}\label{wls3}
	\bar{\theta}(k)=\theta(k)-P^{\frac{1}{2}}(k)\beta(k).
\end{equation}

Finally, the continuous-time estimates used for the design of adaptive control can be defined piecewise as follows:
\begin{equation}\label{wls4}
	\hat{\theta}(t)=\bar{\theta}(k), \text{ for any } t \in (k, k +1],
\end{equation}
for all $k\in\mathbb{N}$. For the above regularized estimates, we have the following lemma, and its proof is given in Appendix.
\begin{lemma}\label{le3}
	If Assumption \ref{ass1} is satisfied, then for any admissible strategy $u(t)$, the family of regularized WLS estimates $\{\hat{\theta}(t),t\ge0\}$ defined by (\ref{wls1})-(\ref{wls4}) has the following properties.
	\par 
	1) Self-convergence, that is, $\hat{\theta}(t)$ converges a.s. to a finite random matrix as $t\to\infty$.
	\par 
	2) The family of matrices $\{A(t),B(t);t\ge0\}$ is uniformly stabilizable and $\{A(t),Q^{1/2};t\ge0\}$ is uniformly detectable where $$[A(t),B(t)]=\hat{\theta}^{\mathsf{T}}(t).$$
	
	3) 
	$\int_{0}^{t}|(\hat{\theta}(s)-\theta)^{\mathsf{T}}\varphi(s)|^2ds=o(r(t))+O(1),$ as $t\to \infty$, 
	where $r(t)$ is defined by (\ref{r}).
	
	4) Boundedness, that is, $\sup\limits_{k\in\mathbb{N}}\|P^{-\frac{1}{2}}(k)(\hat{\theta}(k+1)-\theta)\|<\infty$.
\end{lemma}
\section{Main Results}\label{section3}
For simplicity, we rewrite the estimates given by (\ref{wls4}) as
$$\hat{\theta}^{\mathsf{T}}(t)=[A(t),B(t)].$$

For any $t\ge0$, since $(A(t),B(t))$ is stabilizable and $(A(t),Q^{1/2})$ is detectable,  the following ARE:
\begin{equation}\label{are2}
	A^{\mathsf{T}}(t)X(t)+X(t)A(t)+Q-X(t)S(t)X(t)=0
\end{equation}
has a stabilizing solution $X(t)$, where $S(t)=B(t)R^{-1}B^{\mathsf{T}}(t)$.

By certainty equivalence principle, the adaptive control could have the following form:

\begin{equation}\label{key4}
	u(t)=L(t)x(t),\text{ with } L(t)=-R^{-1}B^{\mathsf{T}}(t)X(t).
\end{equation}

According to  the well-known Fel'dbaum dual principle in optimal control, achieving a balance between control and estimation is crucial for uncertain systems.  A comparable concept in reinforcement learning is the trade-off between ``exploitation and exploration". Building upon these ideas and following the technique in \cite{b2}, some diminishing excitation or exploration signals that are helpful for the strong consistency of the estimates but not essentially influence the performance, are incorporated into the adaptive control, i.e., for $t\in(k,k+1],$
$$
u^*(t)=L(t)x(t)+\gamma_k(v(t)-v(k)) \text{  or}
$$
\begin{equation}\label{c1}
	\mathsf{d} u^*(t)=L(t)\mathsf{d}x(t)+\gamma_k\mathsf{d}v(t),
\end{equation}
where $\gamma_k=(1/k)^{\frac{1}{5}} \text{ with } \gamma_0=0$, and the processes $\{v(t),t \ge 0\}$  are chosen as a sequence of independent standard Wiener processes that are independent of $\{w(t),t\ge0\}$ and $\{\eta_k,k\in\mathbb{N}\}$.

\begin{theorem}\label{th1}{\rm (Global Stability)}
	Let Assumption \ref{ass1} be satisfied. Then, under the
	adaptive control (\ref{c1}), the system (\ref{s1}) is globally stable in the sense that for any initial state $x(0)$,
	\begin{equation}
		\limsup\limits_{T\to \infty}\frac{1}{T}\int_{0}^{T}|x(t)|^2\mathsf{d}t<\infty \quad a.s.
	\end{equation}
\end{theorem}

Next, we discuss whether the estimates $\hat{\theta}(t)$ given by (\ref{wls4}) converge to the true value $\theta$. It should be noted that, under the conditions that $(A,B)$ is controllable and $(A,Q^{1/2})$ is observable, the strong consistency of $\hat{\theta}(t)$ (i.e., $\lim_{t\to\infty}\hat{\theta}(t)=\theta$) was established in \cite{a2}. However, such  strong consistency property may not be  true since we only assume that  $(A,B)$ is stabilizable and $(A,Q^{1/2})$ is detectable in the current paper. 

Let us give some explanations. First, we denote $B_0=[B,D]$ and $n_1$ is the dimension of $\text{Im}(A|B_0)$. Let us choose $n_1$ linear independent bases $e_1,e_2,\cdots,e_{n_1}\in\text{Im}(A|B_0)$, and $n-n_1$ linear independent bases $e_{n_1+1},\cdots,e_{n}\in \text{Im}(A|B_0)^\bot$, and denote $U=[e_1,e_2,\cdots,e_n]$, then the system (\ref{s1}) can be rewritten into
\begin{equation}\label{s2}
	\mathsf{d}\bar{x}=(\bar{A}\bar{x}+\bar{B}u)\mathsf{d}t+\bar{D}\mathsf{d}w
\end{equation}
where $$\bar{x}=U^{-1}x=\left(\begin{array}{cc}
	x_1  \\ 
	x_2 \\ 
\end{array}\right),\bar{A}=U^{-1}AU=\left(\begin{array}{cc}
	A_1 & A_2  \\ 
	0 & A_3\\ 
\end{array}\right),$$
and
$$\bar{B}=U^{-1}B=\left(\begin{array}{cc}
	B_1  \\ 
	0 \\ 
\end{array}\right),\bar{D}=U^{-1}D=\left(\begin{array}{cc}
	D_1  \\ 
	0 \\ 
\end{array}\right).$$
According to the special form of $\bar{A},\bar{B},\bar{D}$,  the dynamics of  $x_2$ is given by 
$$\mathsf{d}x_2=A_3x_2\mathsf{d}t.$$
It is obvious that there exists some initial state $x_2(0)$ (e.g., $x_2(0)=0$) such that we cannot estimate the true value $A_3$ no matter what the input control is. In fact, we could estimate the value $$\left(\left(\begin{array}{cc}
	A_1   \\ 
	0 \\ 
\end{array}\right),\left(\begin{array}{cc}
	B_1  \\ 
	0\\ 
\end{array}\right)\right)$$
for all initial state, which will be shown in the following theorem. But, since $A,B,D$ are unknown, one direct challenge is how to obtain a desired matrix $U$. To address this issue, let us denote a linear space
$$\mathcal{U}(k)=\Big\{\text{span}(\eta)\Big|\eta\text{ is an eigenvector of }\int_{0}^{k}\varphi(s)\varphi^{\mathsf{T}}(s)ds$$$$ \quad\quad\quad\text{whose corresponding eigenvalue } \lambda<\log k\Big\}.$$
Then we have the following lemma and its proof is given in Appendix.
\begin{lemma}\label{le8}
	For the system (\ref{s1}) with the adaptive control (\ref{c1}), there exists some $k_0\ge0$ such that when $k\ge k_0$,  $[\eta_1^{\mathsf{T}},\eta_2^{\mathsf{T}}]^{\mathsf{T}}\in \mathcal{U}(k)$ if and only if  $\eta_1\in \text{Im}(A|B_0)^\bot$ and $\eta_2=0$.
\end{lemma}
From the above lemma, when $k$ is large enough, a space $\text{Im}(A|B_0)^\bot$ can be obtained easily from the space $\mathcal{U}(k)$. Hence, a desired matrix $U$ can be chosen. Then for such $U$, we have the following theorem.
\begin{theorem}\label{th2}{\rm (Strong Consistency)}
	Let Assumption \ref{ass1} be satisfied. Then, for the system (\ref{s1}) with the
	adaptive control (\ref{c1}),  some components of the parameter estimates $\hat{\theta}(t)$ are strongly consistent, i.e.,
	\begin{align}
		&\lim\limits_{t\to\infty}U^{-1}\hat{\theta}^{\mathsf{T}}(t) \mathsf{diag}(U,I_m)\mathsf{diag}(I_{n_1},0_{n-n_1},I_m)\notag\\
		&=U^{-1}\theta^{\mathsf{T}}\mathsf{diag}(U,I_m)\mathsf{diag}(I_{n_1},0_{n-n_1},I_m)\notag\\
		&=\left(\begin{array}{ccc}
			A_1 & 0 & B_1  \\ 
			0 & 0 & 0  \\ 
		\end{array}\right),\quad a.s.
	\end{align}
	where $\theta$ is the true system parameter.
\end{theorem}
\begin{remark}
	If $(A,B)$ is controllable, then the above theorem implies the strong consistency of $\hat{\theta}(t)$, which have been established in Theorem 2 of \cite{a2}. However, when $(A,B)$ is stabilizable, the estimates of $(A_2,A_3)$ are not guaranteed to be strongly consistent. Fortunately, the above property for estimated models is enough for the design of optimal adaptive control, which will be shown in the following theorem.
\end{remark}

\begin{remark}
	From the proof of Theorem \ref{th2}, one can see that as long as the excitation corresponding to the estimated parameters is rich, we have the ability to estimate their true values. It is worth mentioning that the noise is helpful for parameter estimation in a certain sense.
\end{remark}

Note that for $\bar{A},\bar{B}$, the corresponding ARE (\ref{are1}) is equivalent to the following three equations:
\begin{equation}\label{are3}
	A_1^{\mathsf{T}}X_1+X_1A_1+Q_1-X_1B_1R^{-1}B_1^{\mathsf{T}}X_1=0,
\end{equation}
$$(A_1^{\mathsf{T}}-X_1B_1R^{-1}B_1^{\mathsf{T}})X_2+X_2A_3+Q_2+X_1A_2=0,$$
$$A_3^{\mathsf{T}}X_3\!+\!X_3A_3\!+\!X_2^{\mathsf{T}}A_2\!+\!A_2^{\mathsf{T}}X_2\!+\!Q_3\!-\!X_2B_1R^{-1}B_1^{\mathsf{T}}X_2\!=\!0,$$
where
$$\bar{X}\!=\!U^{\mathsf{T}}XU\!=\!\left(\begin{array}{cc}
	X_1 & X_2  \\ 
	X_2^{\mathsf{T}} & X_3 \\ 
\end{array}\right), \bar{Q}\!=\!U^{\mathsf{T}}XU\!=\!\left(\begin{array}{cc}
	Q_1 & Q_2  \\ 
	Q_2^{\mathsf{T}} & Q_3 \\ 
\end{array}\right).$$
Since only the value $(A_1,B_1)$ can be estimated, we may only obtain accurate solution of the first equation (\ref{are3}) when $(A,B)$ are unknown. Fortunately, it is sufficient for the design of optimal adaptive control, since the component of state $x_2$ converges to zero exponentially, and does not matter the performance of both the system (\ref{s1}) and the ergodic cost function (\ref{cost}). 
%In fact, we only need the value $X_1$. 
%For the optimal control defined by (\ref{u})
%$$u=-R^{-1}B^{\mathsf{T}}{X}x=-R^{-1}%%(B_1^{\mathsf{T}}X_1x_1+B_1^{\mathsf{T}}X_2x_2),$$
%since the state $x_2$ converges to zero, the value $X_2$ does not matter the performance of the system (\ref{s1}) and the ergodic cost function (\ref{c1}), which induces the following proposition.
%\begin{proposition}
%The control $u=-R^{-1}B_1^{\mathsf{T}}X_1x_1+L_2x_2$ can stabilize the system (\ref{s1}) and is optimal for the ergodic cost function (\ref{c1}) with $J(u^*)=\text{tr}(D^{\mathsf{T}}XD), a.s.$, where $X_1$ is the stabilizing solution to the ARE (\ref{are3}) and $L_2\in\mathbb{R}^{m\times n-n_1}$ is an arbitrary matrix.
%\end{proposition}

The following theorem indicates that  the proposed adaptive control (\ref{c1}) is optimal.

\begin{theorem}\label{th3}{\rm (Optimality)}
	For the system (\ref{s1}) with the cost function (\ref{cost}), if  Assumption \ref{ass1} is satisfied, then the adaptive control (\ref{c1}) is optimal, that is, 
	\begin{equation}
		J(u^*)=\text{tr}(D^{\mathsf{T}}X^*D)\quad a.s.
	\end{equation}
	where $X^*$ is the stabilizing solution to the ARE (\ref{are1}).
\end{theorem}
\section{Proofs}\label{section4}
\subsection{Proof of Theorem \ref{th1}:}
Before the proof, we need the following lemma.
\begin{lemma}\cite{Liu2024}\label{le9}
	For the processes $\{v(t),t \ge 0\}$ and the sequences $\{\gamma_k, k\in\mathbb{N}\}$ appeared in (\ref{c1}), the following properties hold:
	$$
	\limsup\limits_{N\to \infty}\frac{1}{N}\sum_{k=1}^{N}\int_{k}^{k+1}\gamma_k^2|v(t)-v(k)|^2\mathsf{d} t=0\quad a.s.$$
\end{lemma}

By Lemma \ref{le3}, there are random matrices $A(\infty)$ and $B(\infty)$ such that
$$\lim\limits_{t\to \infty}A(t)=A(\infty), \quad\lim\limits_{t\to \infty}B(t)=B(\infty)\quad a.s.$$
and that $(A(\infty),B(\infty))$ is stabilizable. For simplicity of the remaining descriptions, we denote
\begin{equation}\label{l2}
	\Phi(t)=A(t)+B(t)L(t).
\end{equation}
where $L(t)$  is defined in (\ref{c1}).

It is easy to see that $\Phi(t)$ is uniformly stable and convergent. Hence, by Lyapunov equation, there exist some uniformly bounded positive definite matrices $K(t)$ such that
\begin{equation}\label{key7}
	\Phi^{\mathsf{T}}(t)K(t)+K(t)\Phi(t)=-I.
\end{equation}

Next, we proceed to verify that
\begin{equation}\label{l5}
	\sum\limits_{k=0}^{N}|x(k)|^2=O(N)+o(r(N))\quad a.s.
\end{equation}
where $r(t)$ is defined in (\ref{r}). Note that for $t\in(k,k+1]$ and $k\in\mathbb{N}$, under adaptive strategies (\ref{c1}) the system (\ref{s1}) will be
\begin{align}\label{s}
	\mathsf{d}x(t)&=(Ax(t)+Bu^*(t))\mathsf{d}t+D\mathsf{d}w(t)\\
	&=(\Phi(t)x(t)+\delta(t)+\gamma_k(v_1(t)-v_1(k)))\mathsf{d}t+D\mathsf{d}w(t),\notag
\end{align}
where $\delta(t)=(\theta-\hat{\theta}(t))^{\mathsf{T}}\varphi(t)$ and $v_1(t)=Bv(t)$.

Then, it follows that
$$x(k+1)=e^{\Phi(k)}x(k)+\int_{k}^{k+1}e^{(k+1-t)\Phi(k)}Ddw(t)$$
$$+\int_{k}^{k+1}e^{(k+1-t)\Phi(k)}(\delta(t)+\gamma_k(v_1(t)-v_1(k)))dt.$$

Since $\Phi(k)$ is uniformly stable and convergent, by Cauchy-Schwarz inequality, it is easy to get
$$	|x(k+1)|^2\le m|x(k)|^2+m_1(\int_{k}^{k+1}e^{(k+1-t)\Phi(k)}Ddw(t))^2$$
$$+m_2(\int_{k}^{k+1}|\delta(t)|^2dt+\int_{k}^{k+1}\gamma_k^2|v(t)-v(k)|^2dt),$$
where $0<m<1$ and $m_1,m_2>0$ are some fixed constants associated with the supremum of the family $\{e^{\Phi(k)}\}$. Then, it is easy to see that
\begin{align}
	&\frac{1}{N}(1-m)\sum\limits_{k=1}^{N}|x(k)|^2\notag\\
	&=O(\frac{1}{N}\sum\limits_{k=1}^{N}(\int_{k}^{k+1}e^{(k+1-t)\Phi(k)}Ddw(t))^2)+\notag\\
	&O(\frac{1}{N}\int_{1}^{N}|\delta(t)|^2dt)+O(\frac{1}{N}\sum_{k=1}^{N}\int_{k}^{k+1}\gamma_k^2|v_1(t)-v_1(k)|^2dt)\notag\\
	&=O(1)+o(\frac{1}{N}r(N)),\notag
\end{align}
where the first part can use Lemma 1 (Etemadi) of 5.2 in \cite{b7}, the second part is the direct result of Lemma \ref{le3} 3) and the third part can be estimated by the consequence of Lemma \ref{le9}. 

Finally, we proceed to prove Theorem \ref{th1}. Applying the Ito's formula to $\langle K(t)x(t),x(t)\rangle$ where $\langle\cdot,\cdot\rangle$ represents the inner product, and noting that $K(t)$ defined in (\ref{key7}) is actually constant in any interval $t\in(i,i+1]$ for any $i\in\mathbb{N}$, it follows that
$$\mathsf{d}\langle K(t)x(t),x(t)\rangle=\text{tr}(K(t)DD^{\mathsf{T}})\mathsf{d}t+2\langle K(t)x(t),D\mathsf{d}w(t)\rangle$$
$$
+2\langle K(t)x(t),\Phi(t)x(t)+\delta(t)+\gamma_i(v_1(t)-v_1(i))\rangle \mathsf{d}t,
$$
which in conjunction with equation (\ref{key7}) gives
$$\mathsf{d}\langle K(t)x(t),x(t)\rangle+|x(t)|^2\mathsf{d}t=\text{tr}(K(t)DD^{\mathsf{T}})\mathsf{d}t+$$
\begin{equation}\label{e2}
	2\langle K(t)x(t),D\mathsf{d}w(t)\rangle+2\langle K(t)x(t),\delta(t)+\gamma_i(v_1(t)-v_1(i))\rangle \mathsf{d}t.
\end{equation}

For the second part on the right-hand-side, by the boundedness of $K(t)$, it follows from Lemma 12.3 of \cite{b2} that
\begin{equation}\label{z1}
	|\int_{0}^{t}\langle K(t)x(t),Ddw(t)\rangle|=O\big([\int_{0}^{t}|x(t)|^2dt]^{\frac{1}{2}+\epsilon}\big),
\end{equation}
for any $\epsilon\in(0,1/2)$.

By Lemma \ref{le3}, it follows that
\begin{equation}\label{z2}
	\int_{0}^{t}|\delta(s)|^2ds=o(r(t))+O(1).
\end{equation}

Then, integrating the equation (\ref{e2}) over the interval $(0,T)$, and using (\ref{z1})-(\ref{z2}), Lemma \ref{le9} and the Cauchy-Schwarz inequality, it follows that
$$\sum\limits_{i=0}^{\lfloor T\rfloor-1}\big(\langle K(i)x(i+1),x(i+1)\rangle-\langle K(i)x(i),x(i)\rangle\big)$$
$$+\langle K(\lfloor T\rfloor)x(T),x(T)\rangle-\langle K(\lfloor T\rfloor)x(\lfloor T\rfloor),x(\lfloor T\rfloor)\rangle$$
$$+\int_{0}^{T}|x(t)|^2dt= \int_{0}^{T}\text{tr}(K(t)DD^{\mathsf{T}})dt+$$
$$O\big((\int_{0}^{T}|x(t)|^2dt)^{\frac{1}{2}+\epsilon}\big)+(\int_{0}^{T}|x(t)|^2dt)^\frac{1}{2}\cdot o(T^{\frac{1}{2}})+$$
\begin{equation}\label{a12}
	(\int_{0}^{T}|x(t)|^2dt)^\frac{1}{2}(o(r(T))+O(1))^\frac{1}{2},
\end{equation}
where $\lfloor T\rfloor$ denotes the integer part of $T$.

Since $K(t)$ is uniformly bounded, we have
$$\sum\limits_{i=0}^{\lfloor T\rfloor-1}\big(\langle K(i)x(i+1),x(i+1)\rangle-\langle K(i)x(i),x(i)\rangle\big)$$
$$+\langle K(\lfloor T\rfloor)x(T),x(T)\rangle-\langle K(\lfloor T\rfloor)x(\lfloor T\rfloor),x(\lfloor T\rfloor)\rangle$$
\begin{equation}\label{ee3}
	=O(\sum\limits_{i=0}^{\lfloor T\rfloor}|x(i)|^2)+\langle K(\lfloor T\rfloor)x(T),x(T)\rangle.
\end{equation}

By Lemma \ref{le9} and the convergence of $L(t)$, it follows that
\begin{equation}\label{ee4}
	r(T)=\|P^{-1}(0)\|+\int_{0}^{T}|\varphi(s)|^2ds=O(\int_{0}^{T}|x(t)|^2dt).
\end{equation}

Finally, by (\ref{l5}) and (\ref{ee3})-(\ref{ee4}), the equality (\ref{a12}) is simplified to
$$\langle K(\lfloor T\rfloor)x(T),x(T)\rangle+\int_{0}^{T}|x(t)|^2dt$$
\begin{equation}\label{a13}
	=O(T)+o(\int_{0}^{T}|x(t)|^2dt)+\int_{0}^{T}\text{tr}(K(t)DD^{\mathsf{T}})dt,
\end{equation}
which implies the desired result of Theorem \ref{th1}. 
\subsection{Proof of Theorem \ref{th2}:}
Note that the augmented state 
$$\varphi(t)=[x^{\mathsf{T}}(t),u^{\mathsf{T}}(t)]^{\mathsf{T}}$$ 
satisfies the following equation for any $t\in(k,k+1]$ and any $k\in\mathbb{N}$:
\begin{equation}\label{a15}
	\mathsf{d} \varphi(t)=C_k\varphi(t)\mathsf{d} t+G_k\mathsf{d}\xi(t),
\end{equation}
where
$$C_k=\left(\begin{array}{cc}
	A & B \\ 
	L(k)A & L(k)B  \\ 
\end{array}\right),$$
$$G_k=\left(\begin{array}{cc}
	D & 0  \\ 
	L(k)D & \gamma_kI_{m}  \\ 
\end{array}\right),$$
$$\xi(t)=(w(t),v(t)),$$

Now, we need the following lemmas, and some proofs are given in Appendix.

\begin{lemma}\cite{a2}\label{le4}
	Let $\{H(k)-H(t),t\in(k,k+1],k\in\mathbb{N}\}$ be given by
	$$H(t)-H(k)=\int_{k}^{t}\varphi(s)\varphi^{\mathsf{T}}(s)\mathsf{d} s,$$then
	\begin{align}
		&H(k+1)-H(k)\ge\int_{k}^{k+1}e^{C_k(k+1-s)}[G_kG_k^{\mathsf{T}}\int_{k}^{s}\mathsf{d}\tau\notag\\
		&\quad\quad\quad\quad\quad\quad\quad\quad\quad+\int_{k}^{s}\mathsf{d}g(\tau)]e^{C^{\mathsf{T}}_k(k+1-s)}\mathsf{d} s\notag\
	\end{align}
	where
	$$g(s)=\varphi(s)\mathsf{d}\xi^{\mathsf{T}}(s)G_k^{\mathsf{T}}+G_k\mathsf{d}\xi(s)\varphi^{\mathsf{T}}\text{ for }s\in(k,k+1].$$
\end{lemma}
\begin{lemma}\label{le5}
	The process $\{M(k+1)\!-\!M(k),k\in\mathbb{N}\}$ is given by
	\begin{align}
		&M(k+1)-M(k)\notag\\
		&=\int_{k}^{k+1}e^{C_k(k+1-s)}\int_{k}^{s}\mathsf{d}g(\tau)e^{C^{\mathsf{T}}_k(k+1-s)}\mathsf{d} s,\notag
	\end{align}
	then, for any sufficiently large $k\in\mathbb{N}$, the following property holds for any $\epsilon\in(0,1/2)$,
	$$\lambda_{\max}(M(k+1))=O(\int_{0}^{k}|x(s)|^2\mathsf{d}s)^{1/2+\epsilon}\quad a.s.$$
\end{lemma}
\begin{lemma}\label{le6}
	For the system (\ref{s1}) with the adaptive control (\ref{c1}), we have for all sufficiently large $k\in\mathbb{N}$,
	$$\lambda_{\min}\left(\int_{0}^{k}U_1\varphi(s)\varphi^{\mathsf{T}}(s)U_1^{\mathsf{T}}\mathsf{d} s\right)>k^{1/2}\quad a.s.$$
	$$\lambda_{\max}\left(\int_{0}^{k}U_2\varphi(s)\varphi^{\mathsf{T}}(s)U_2^{\mathsf{T}}\mathsf{d} s\right)=O(1)\quad a.s.$$
	where 
	$$U_1=\left(\begin{array}{ccc}
		I_{n_1} & 0 &0  \\ 
		0& 0  &I_m\\ 
	\end{array}\right)\mathsf{diag}(U^{-1},I_m),$$$$U_2=[0,I_{n-n_1},0]\mathsf{diag}(U^{-1},I_m).$$
\end{lemma}

\begin{lemma}\label{le7}
	For the WLS estimator (\ref{wls1})-(\ref{wls4}), let us denote 
	$$\theta=\left(\begin{array}{cc}
		\theta_1  \\ 
		\theta_2 \\ 
	\end{array}\right),\hat{\theta}(t)=\left(\begin{array}{cc}
		\theta_1(t)  \\ 
		\theta_2(t) \\ 
	\end{array}\right),\varphi(t)=\left(\begin{array}{cc}
		\varphi_1(t)  \\ 
		\varphi_2(t) \\ 
	\end{array}\right)$$
	if for any sufficiently large $k\in\mathbb{N}$ we have 
	$$\frac{\lambda_{\max}(I+\int_{0}^{k}\varphi_2(s)\varphi_2^{\mathsf{T}}(s)ds)}{a(t)\lambda_{\min}(\int_{0}^{k}\varphi_1(s)\varphi_1^{\mathsf{T}}(s)ds)}=o(1)\quad a.s.$$
	then the follow property holds:
	$$\lim\limits_{t\to\infty}\theta_1(t)=\theta_1\quad a.s.$$
\end{lemma}
Noting that $$1/a(t)=O(\log^2r(t))=O(\log^2 t),$$ 
and from Lemma \ref{le6}, applying Lemma \ref{le7} to $\hat{\theta}^{\mathsf{T}}(t) \mathsf{diag}(U,I_m)$,  
it follows that
\begin{align}
	&\lim\limits_{t\to\infty}\hat{\theta}^{\mathsf{T}}(t) \mathsf{diag}(U,I_m)\mathsf{diag}(I_{n_1},0_{n-n_1},I_m)\notag\\
	&=\theta^{\mathsf{T}} \mathsf{diag}(U,I_n)\mathsf{diag}(I_{n_1},0_{n-n_1},I_m)\quad a.s.\notag\\
	&=U\left(\begin{array}{ccc}
		A_1 & A_2 & B_1  \\ 
		0 & A_3 & 0  \\ 
	\end{array}\right)\mathsf{diag}(I_{n_1},0_{n-n_1},I_m),\notag\\
	&=U\left(\begin{array}{ccc}
		A_1 & 0 & B_1  \\ 
		0 & 0 & 0  \\ 
	\end{array}\right),
\end{align}
which completes the proof.

\subsection{Proof of Theorem \ref{th3}:}
Before the proof, we need the following lemma.
\begin{lemma}\cite{Liu2024}\label{le10}
	Assume that
	\begin{equation}\notag
		\limsup\limits_{T\to \infty}\frac{1}{T}\int_{0}^{T}|x(t)|^2\mathsf{d} t<\infty \quad  a.s.
	\end{equation}
	and 
	\begin{equation}\notag
		\lim\limits_{t\to \infty}V(t)=0 \quad  a.s.
	\end{equation}
	then the following equality is true:
	\begin{equation}\notag
		\lim\limits_{t\to \infty}\frac{1}{T}\int_{0}^{T}\langle V(t)x(t),x(t)\rangle \mathsf{d} t=0 \quad  a.s.
	\end{equation}
\end{lemma}
By Theorem \ref{th1} and Lemma \ref{le9}, it is easy to verify that the adaptive control (\ref{c1}) is admissible. Next, we proceed to verify that
\begin{equation}
	J(u^*)=\text{tr}(D^{\mathsf{T}}X^*D).
\end{equation}
Under the adaptive control (\ref{c1}), the system (\ref{s2}) will be
\begin{align}\label{s3}
	\mathsf{d}x_1&=(A_1x_1+A_2x_2+B_1u^*(t))\mathsf{d}t+D_1\mathsf{d}w(t)\notag\\
	&=(\Phi_1(t)x_1+\Phi_2(t)x_2\notag\\
	&\quad\quad+\gamma_{ \lfloor t\rfloor }(v_1(t)-v_1( \lfloor t\rfloor )))\mathsf{d} t+D_1\mathsf{d} w(t),
\end{align}
with $\mathsf{d}x_2=A_3x_2\mathsf{d}t$, where $\Phi_1(t)=A_1-B_1L_1(t)$, $\Phi_2(t)=A_2-B_1L_2(t)$ with $[L_1(t),L_2(t)]=L(t)$ and $v_1(t)=B_1v(t)$.

Applying Ito's formula to $\langle X_1^*x_1(t),x_1(t)\rangle$, where $X_1^*$ is the stabilizing solution to the ARE (\ref{are3}), it follows that
\begin{align}\label{e1}
	&\mathsf{d}\langle X_1^*x_1(t),x_1(t)\rangle\notag\\
	&=2\langle X_1^*x_1(t),\Phi_1(t)x_1\!+\!\Phi_2(t)x_2\!+\!\gamma_{ \lfloor t\rfloor }(v_1(t)\!-\!v_1( \lfloor t\rfloor ))\rangle \mathsf{d} t\notag\\
	&\quad+\text{tr}(X_1^*D_1D_1^{\mathsf{T}})\mathsf{d} t+2\langle X_1^*x_1(t),D_1\mathsf{d} w(t)\rangle.
\end{align}

By integrating (\ref{e1}) over the interval $(0,T)$, we have
\begin{align}\label{z3}
	&\langle X_1^*x_1(T),x_1(T)\rangle-\langle X_1^*x_1(0),x_1(0)\rangle\notag\\
	&=2\int_{0}^{T}\langle X_1^*x_1(t),\Phi_1(t)x_1\!+\!\Phi_2(t)x_2\rangle \mathsf{d} t\notag\\
	&\quad\quad+2\int_{0}^{T}\langle X_i^*x(t),\gamma_{ \lfloor t\rfloor }(v_1(t)\!-\!v_1( \lfloor t\rfloor ))\rangle \mathsf{d} t\notag\\
	&\quad\quad+T\text{tr}(X_1^*DD^{\mathsf{T}})+2\int_{0}^{T}\langle X_i^*x(t),\mathsf{d} w(t)\rangle.
\end{align}
Now, we analyze the right-hand side of the above equation. By the Cauchy-Schwarz inequality and Lemma \ref{le9}, it follows that
\begin{align}\label{k4}
	&\int_{0}^{T}\langle X_1^*x_1(t),\Phi_2(t)x_2\rangle \mathsf{d} t\notag\\
	&\le(\int_{0}^{T}|X_1^*x_1(t)|^2\mathsf{d} t)^\frac{1}{2}(\int_{0}^{T}|\Phi_2(t)x_2|^2\mathsf{d} t)^\frac{1}{2}\notag\\
	&=O\Big((\int_{0}^{T}|x(t)|^2\mathsf{d} t)^\frac{1}{2}\Big)
\end{align}
and
\begin{align}\label{k5}
	&\int_{0}^{T}\langle X_1^*x_1(t),\gamma_{ \lfloor t\rfloor }(v_1(t)-v_1( \lfloor t\rfloor ))\rangle \mathsf{d} t\notag\\
	&\le(\int_{0}^{T}|X_1^*x_1(t)|^2\mathsf{d} t)^\frac{1}{2}(\int_{0}^{T}\gamma_{ \lfloor t\rfloor }^2|v_1(t)\!-\!v_1( \lfloor t\rfloor )|^2\mathsf{d} t)^\frac{1}{2}\notag\\
	&=O\big((\int_{0}^{T}|x(t)|^2\mathsf{d} t)^\frac{1}{2}\big)\cdot o(T^\frac{1}{2})
\end{align}
By  Lemma 12.3 of \cite{b2}, it follows that for any $\epsilon\in(0,1/2),$
\begin{equation}\label{k6}
	|\int_{0}^{t}\langle X_1^*x_1(t),D\mathsf{d} w(t)\rangle|=O\Big((\int_{0}^{t}|x_1(t)|^2\mathsf{d} t)^{\frac{1}{2}+\epsilon}\Big),
\end{equation}
By (\ref{k4})-(\ref{k6}), (\ref{z3}) is simplified to
\begin{equation}\label{ll6}
	\limsup\limits_{T\to \infty}\frac{1}{T}\int_{0}^{T}\langle -2X_1^*x_1(t),\Phi_1(t)x_1(t)\rangle \mathsf{d} t=\text{tr}(X_1^*D_1D_1^{\mathsf{T}}).
\end{equation}
By Lemma \ref{le9}, Lemma \ref{le10} and the CARE (\ref{are3}), it follows that
\begin{align}
	&J(u^*(t))\notag\\
	&=\limsup\limits_{T\to \infty}\frac{1}{T}\int_{0}^{T}\big(\langle Qx(t),x(t)\rangle\!+\!\langle Ru^*(t),u^*(t)\rangle\mathsf{d} t\notag\\
	&=\limsup\limits_{T\to \infty}\frac{1}{T}\int_{0}^{T}\big(\langle Q_1x_1(t),x_1(t)\rangle\!+\!\langle Ru^*(t),u^*(t)\rangle\mathsf{d} t\notag\\
	&=\limsup\limits_{T\to \infty}\frac{1}{T}\int_{0}^{T}\big(\langle (Q_1+X_1^*B_1R^{-1}B_1^{\mathsf{T}}X_1^*)x_1(t),x_1(t)\rangle\mathsf{d} t\notag\\
	&=\limsup\limits_{T\to \infty}\frac{1}{T}\int_{0}^{T}\langle -2X_1^*x_1(t),(A-B_1R^{-1}B_1^{\mathsf{T}}X_1^*)x_1(t)\rangle \mathsf{d} t \notag\\
	&=\limsup\limits_{T\to \infty}\frac{1}{T}\int_{0}^{T}\langle -2X_1^*x_1(t),\Phi_1(t)x_1(t)\rangle \mathsf{d} t \notag\\
	&=\text{tr}(X_1^*D_1D_1^{\mathsf{T}})=\text{tr}(X^*DD^{\mathsf{T}}).
\end{align}
Hence, Theorem \ref{th3} is true.

\section{Conclusions}\label{section6}
The continuous-time adaptive LQG control problem is one of the well-known problems in adaptive control. In this paper, a complete solution is provided for this problem under the assumptions that $(A, B)$ is stabilizable and $(A, Q^{1/2})$ is detectable. The major technical methods associated with this problem are: a
WLS algorithm with random regularization method ensures that a family of estimated models for the unknown system matrices are self-convergent, uniformly stabilizable and detectable respectively; a certainty equivalence control is used to ensure the performance of the closed-loop adaptive control system; and a suitable diminishing excitation is incorporated in the adaptive control to ensure strong consistency of the desired components of the estimates. Even if not all components of the estimated models converge to the true values, the proposed adaptive control still has the  ability to achieve global stability and optimality. It is worth mentioning that, for the adaptive discrete-time LQ control problem, the same results as those established in this paper can be obtained similarly. However, many fascinating issues remain to be further investigated. For example, how to solve the adaptive problem where the unknown system parameters are time-varying? 

\section*{Appendix}\label{appendix}
\textbf{Proof of Lemma \ref{le2}}:
It is obvious that $1)\iff2)$ and $3)\iff4)$. Hence, we only need to prove that 2)  $(sI-A,B)$ is full of row rank for any $s\in\sigma^+(A)$
if and only if 3) $Z(s,A,B)$ is full of row rank for any $s\in\sigma^+(A)$. For this purpose, note that
$$[\eta_1^\mathsf{T},\eta_2^\mathsf{T}]Z(s,A,B)=0,$$
if and only if
$$(\eta_1^\mathsf{T}+\eta_2^\mathsf{T}j)(sI-A,B)=0,$$
since they are equivalent to
\begin{align}\notag
	\left\{
	\begin{aligned}
		&s_1\eta_1^\mathsf{T}-s_2\eta_2^\mathsf{T}-\eta_1^\mathsf{T}A=0,\\
		&s_1\eta_2^\mathsf{T}+s_2\eta_1^\mathsf{T}-\eta_2^\mathsf{T}A=0,\\
		&\eta_1^\mathsf{T}B=\eta_2^\mathsf{T}B=0,
	\end{aligned}
	\right.
\end{align}
where $s=s_1+s_2j$ and $j$ is the imaginary unit. Hence, the proof is completed.

\vspace{5mm}
\textbf{Proof of Lemma \ref{le3}}:
To vefify 1), since $\theta(t)$ and $P(t)$ are convergent, we only need to verify that $\{\beta(k),k\in\mathbb{N}\}$ is convergent, which will be shown in the proof of 2). 

To verify 2), it suffices from (\ref{wls4}) to show that the sequence $\{\bar{A}(k),\bar{B}(k),Q^{1/2};k\in\mathbb{N}\}$ is uniformly stabilizable and detectable respectively. We proceed to prove it as the following steps, and a similar proof could be found in \cite{a2}.

Step 1: We first prove that for any $\beta(k)\in\mathbb{R}^{(m+n)\times n}$, the pair $(\bar{A}(k),\bar{B}(k))$ is stabilizable and the pair $(\bar{A}(k),Q^{1/2})$ is detectable (or $f_k(\beta(k))\neq0$) almost surely with respect to $\beta(k)$. Since the uncontrollable set is a zero measure set (see, e.g., \cite{b6}), it is trivial.

Step 2: Next, we prove that there exists a positive random variable $\delta_{\infty}$ such that
$$\limsup\limits_{k\to\infty}f_k(\eta_k)\ge \delta_{\infty}\quad a.s.$$
Let us denote
$$\delta_k\triangleq\max_{x\in \mathbb{D}}f_k(x),$$
and
$$\mathbb{D}_k\triangleq\{x\in\mathbb{D}:f_k(x)>\frac{\delta_k}{2}\},$$
where $\mathbb{D}$ the unit ball of $\mathbb{R}^{(n+m)\times n}$ with Euclidean norm.

Note that $\theta(k),P(k),\delta_k$ are all $\mathcal{G}_k\triangleq\sigma(w(s),\eta(s-1);s\le k)$ measurable. Then, we have 
\begin{align}
	P\left(f_k(\eta_k)\ge\frac{\delta_k}{2}|\mathcal{G}_k\right)&=\int_{x\in\mathbb{D}}I(f_k(\eta_k\ge\frac{\delta_k}{2})\mu(dx)\\
	&=\int_{x\in\mathbb{D}_k}\mu(dx)\\
	&=\mu(\mathbb{D}_k),
\end{align}
where $\mu(\cdot)$ is a uniform probability measure defined on $\mathbb{D}$.

Since $\mu(\cdot)$ is absolutely continuous with respect to Lebesgue measure $m(\cdot)$ on $\mathbb{R}^{(m+n)\times n}$, it follows that
$$\mu(x\in D:f_k(x)=0)=0\quad a.s.$$

Note that both $\{\theta(k)\}$ and $\{P(k)\}$ are convergent a.s., we can define a function $f(x)$ as
$$f(x)=\lim\limits_{k\to\infty}f_k(x)\quad a.s.\quad \text{ for } x\in\mathbb{R}^{(n+m)\times n}.$$ 

%Since $(A,B)$ is stabilizable and $(A,Q^{1/2})$ is detectable, it follows that $f(\beta^*)>0$, where $\beta^*=\lim_{t\to\infty}\beta^*(t)$ defined in (\ref{b}). 
Furthermore, it is easy to see that $f_k(x)$ converges  to $f(x)$ uniformly on $\mathbb{D}$. Then it follows that  $f(x)\neq0$ almost surely.  Consequently, $\delta_k\to\delta_\infty\triangleq\max_{x\in D}f(x)>0.$ 

Since $f(x)$ is a continuous function, we have $m(\mathbb{D}_{\infty})>0$, where $\mathbb{D}_{\infty}$ is defined by $$\mathbb{D}_{\infty}=\{x\in\mathbb{D}:f(x)\ge\lambda\delta_{\infty}\} \text{ with }1/2<\lambda<1,$$
which implies that for all sufficient large $k$, $\mu(\mathbb{D}_k)\neq0$, since $\mu(\mathbb{D}_k)=\frac{m(\mathbb{D}_k)}{m(\mathbb{D})}.$

Hence, we know that
$$\sum_{t=1}^{\infty}P\left(f_k(\eta_k)\ge\frac{\delta_k}{2}|\mathcal{G}_k\right)=\infty\quad a.s.$$
By the Borel-Cantelli-Levy Lemma, it follows that
$$\sum_{t=1}^{\infty}I\left(f_k(\eta_k)\ge\frac{\delta_k}{2}\right)=\infty\quad a.s.$$
which implies that
$$\limsup\limits_{k\to\infty}f_k(\eta_k)\ge\frac{1}{2}\lim\limits_{k\to\infty}\delta_k=\frac{1}{2}\delta_{\infty}>0\quad a.s.$$

Step 3: We now prove that there exists positive random variable $\delta$ and $k_0$ such that
$$f(\beta(k))\ge\delta\quad a.s.\quad \forall k\ge k_0$$ 

It is obvious that 
$$\limsup\limits_{k\to\infty}f_k(\beta(k))\ge\frac{\delta_{\infty}}{\gamma}\quad a.s.$$

It is easy to see that there exists a positive random variable $\delta>0$ and a random time $k_0$ such that
$$f_{k_0}(\beta(k_0))\ge 2\delta>0$$
and for any $ k\ge k_0$,
$$|f(\beta(s))-f_k(\beta(s))|\le(\gamma-1)^2\delta\text{ for any }s\in\mathbb{N}.$$
By induction, it is easy to see that
$$f(\beta(k))\ge\delta\quad a.s.\quad \text{ for any } k\ge k_0.$$ 
If  $\beta(k+1)=\beta(k)$, it is trivial; otherwise, $f_{k+1}(\beta(k+1))\ge\gamma f_{k}(\beta(k))$ and it follows that
\begin{align}
	&f(\beta(k+1))\notag\\
	&=f(\beta(k+1))-f_{k+1}(\beta(k+1))+f_{k+1}(\beta(k+1))\notag\\
	&\ge -(\gamma-1)^2\delta+\gamma(\delta-(\gamma-1)^2\delta)\notag\\
	&\ge \delta.\notag
\end{align}

Step 4: Let us first show that $\{\beta(k)\}$ is convergent a.s. It is true since $f_k(x)$ is bounded on $\mathbb{D}$. Then, we have $\lim\limits_{t\to\infty}f(\beta(k))=f>0, a.s.$

Hence, Lemma \ref{le3} 2) is completed.

The property of 3) is the same as the property in Lemma 2 of \cite{a2}, so we omit it. To verify 4), note that 
\begin{align}
	&\|P^{-\frac{1}{2}}(k)(\hat{\theta}(k+1)-\theta)\|\notag\\
	&=\|P^{-\frac{1}{2}}(k)(\theta(k)-\theta-P^{\frac{1}{2}}(k)\beta(k))\|\notag\\
	&\le\|P^{-\frac{1}{2}}(k)(\theta(k)-\theta)\|+\|\beta(k)\|=O(1),\notag
\end{align}
where we have used Lemma \ref{le1}. Hence, the proof is completed.

\vspace{5mm}
\textbf{Proof of Lemma \ref{le8}}: First we consider ``if" part.
If $\eta\in \text{Im}(A|B)^{\perp}$ with $\|\eta\|=1$, then it is easy to verify that 
\begin{align}
	&[\eta^{\mathsf{T}},0]\int_{0}^{k}\varphi(s)\varphi^{\mathsf{T}}(s)\mathsf{d}s[\eta^{\mathsf{T}},0]^{\mathsf{T}}\notag\\
	&=[\bar{\eta}^{\mathsf{T}},0]\mathsf{diag}(U^{-1},I_m)\int_{0}^{k}\varphi(s)\varphi^{\mathsf{T}}(s)\mathsf{d}s\notag\\
	&\quad\quad\cdot\mathsf{diag}(U^{-1},I_m)^{\mathsf{T}}[\bar{\eta}^{\mathsf{T}},0]^{\mathsf{T}}\notag\\
	&=\int_{0}^{k}\eta_1^{\mathsf{T}}x_2(s)x_2^{\mathsf{T}}(s)\eta_1\mathsf{d}s\notag\\
	&=O(1)<\log k,
\end{align}
where $\bar{\eta}^{\mathsf{T}}=\eta^{\mathsf{T}}U=[0,\eta_1^{\mathsf{T}}]$.

Next, we consider ``only if" part. By Lemma \ref{le6}, it follows that if $\eta\in \text{Im}(A,B_0)\oplus R^m$ with $\eta\neq0$, then we have for all sufficiently large $k$,
$$\eta^{\mathsf{T}}\int_{0}^{k}\varphi(s)\varphi^{\mathsf{T}}(s)\mathsf{d}s\eta>\log k.$$
For any $\eta\in\mathcal{U}(k)$ with $\|\eta\|=1$, we can rewrite it as $\eta=[\eta^{\mathsf{T}}_1,0]^{\mathsf{T}}+\eta_2$ where $\eta_1\in \text{Im}(A|B)^{\perp}$ and $\eta_2\in \text{Im}(A,B_0)\oplus R^m$. If $\eta_2\neq0$, then it follows that for all sufficiently large $k$,
\begin{align}
	\eta^{\mathsf{T}}\int_{0}^{k}\varphi(s)\varphi^{\mathsf{T}}(s)\mathsf{d}s\eta
	&\ge\eta_2^{\mathsf{T}}\int_{0}^{k}\varphi(s)\varphi^{\mathsf{T}}(s)\mathsf{d}s\eta_2\notag\\
	&>\log k,\notag
\end{align}
which is contradictory. Hence, the proof is completed.

\vspace{5mm}
\textbf{Proof of Lemma \ref{le5}}:
Before the proof, let $Vr(A)$ denote expanding the matrix $A$ into a column vector by row and $\otimes$ represents the Kronecker product. Then, notice that
\begin{align}
	&\lambda_{\max}(M(i+1))\\
	&=\max\limits_{|\eta|=1}\sum\limits_{k=0}^{i}\eta^{\mathsf{T}}(M(k+1)\!-\!M(k))\eta\notag\\
	&=\max\limits_{|\eta|=1}\sum\limits_{k=0}^{i}\eta^{\mathsf{T}}(\int_{k}^{k+1}e^{C_k(k+1-s)}\notag\\
	&\quad\quad\quad\quad\quad\quad\quad\cdot\int_{k}^{s}\mathsf{d}g(\tau)e^{C^{\mathsf{T}}_k(k+1-s)}\mathsf{d} s)\eta\notag\\
	&=\max\limits_{|\eta|=1}\sum\limits_{k=0}^{i}Vr (\eta^{\mathsf{T}}(\int_{k}^{k+1}e^{C_k(k+1-s)}\notag\\
	&\quad\quad\quad\quad\quad\quad\quad\cdot\int_{k}^{s}\mathsf{d}g(\tau)e^{C^{\mathsf{T}}_k(k+1-s)}\mathsf{d} s)\eta )\notag\\
	&=\max\limits_{|\eta|=1}\sum\limits_{k=0}^{i}(\int_{k}^{k+1}(\eta^{\mathsf{T}}e^{C_k(k+1-s)})\otimes(e^{C^{\mathsf{T}}_k(k+1-s)}\eta)^{\mathsf{T}}\notag\\
	&\quad\quad\quad\quad\quad\quad\quad\cdot\int_{k}^{s}\mathsf{d}Vr(g(\tau))\mathsf{d} s)\notag\\
	&=\max\limits_{|\eta|=1}\sum\limits_{k=0}^{i}(\int_{k}^{k+1}(\eta^{\mathsf{T}}e^{C_k(k+1-s)})\otimes(e^{C^{\mathsf{T}}_k(k+1-s)}\eta)^{\mathsf{T}}\notag\\
	&\quad\quad\quad\quad\quad\quad\quad\cdot\int_{k}^{s}(\varphi\otimes G_k+G_k\otimes\varphi)\mathsf{d}\xi(\tau)\mathsf{d} s)\notag\\
	&=\max\limits_{|\eta|=1}\sum\limits_{k=0}^{i}(\int_{k}^{k+1}\int_{s}^{k+1}(\eta^{\mathsf{T}}e^{C_k(k+1-s)})\otimes(e^{C^{\mathsf{T}}_k(k+1-s)}\eta)^{\mathsf{T}}\notag\\
	&\quad\quad\quad\quad\quad\quad\quad\cdot \mathsf{d}s(\varphi\otimes G_k+G_k\otimes\varphi)\mathsf{d}\xi(\tau))\notag\\
	&=O(\int_{0}^{i}|\varphi(s)|^2\mathsf{d}s)^{1/2+\epsilon}=O(\int_{0}^{i}|x(s)|^2\mathsf{d}s)^{1/2+\epsilon}
\end{align}
where in the last equation we have used the convergence of $(C_k,G_k)$ and Lemma 12.3 of \cite{b2} for any $\epsilon\in(0,1/2)$.

\vspace{5mm}
\textbf{Proof of Lemma \ref{le6}}:
Let $\{N(k+1)-N(k),k\in\mathbb{N}\}$ be given by
\begin{align}
	&N(k+1)-N(k)\notag\\
	&=\int_{k}^{k+1}e^{C_k(k+1-s)}G_kG_k^{\mathsf{T}}\int_{k}^{s}d\tau\, e^{C^{\mathsf{T}}_k(k+1-s)}\mathsf{d} s,\notag
\end{align}
then, we verify that there exists some $c>0$ such that
\begin{equation}\label{n}
	\liminf_{k\to\infty}\gamma_k^{-2}\lambda_{\min}\big(U_1(N(k+1)-N(k))U_1^{\mathsf{T}}\big)>c.
\end{equation}
Note that
\begin{align}
	&G_kG_k^{\mathsf{T}}\!\notag\\
	&=\gamma_k^2\bar{G}_k\bar{G}_k^{\mathsf{T}}\!+\!(1-\gamma_k^2)\left(\begin{array}{cc}
		D & 0  \\ 
		F(k)D & 0   \\ 
	\end{array}\right)\left(\begin{array}{cc}
		D & 0  \\ 
		F(k)D & 0   \\ 
	\end{array}\right)^{\mathsf{T}}\notag\\
	&\ge\gamma_k^2\bar{G}_k\bar{G}_k^{\mathsf{T}},\notag
\end{align}
where
$$\bar{G}_k=\left(\begin{array}{cc}
	D & 0  \\ 
	F(k)D & I_m   \\ 
\end{array}\right).$$
Then, it follows that
\begin{align}
	&N(k+1)-N(k)\notag\\
	&=\int_{k}^{k+1}e^{C_k(k+1-s)}G_kG_k^{\mathsf{T}}\int_{k}^{s}d\tau\, e^{C^{\mathsf{T}}_k(k+1-s)}\mathsf{d}  s,\notag\\
	&\ge\gamma_k^2\int_{k}^{k+1}e^{C_k(k+1-s)}\bar{G}_k\bar{G}_k^{\mathsf{T}}\int_{k}^{s}d\tau\, e^{C^{\mathsf{T}}_k(k+1-s)}\mathsf{d}  s,\notag\\
	&\ge\gamma_k^2\int_{0}^{1}(1-t)e^{C_kt}\bar{G}_k\bar{G}_k^{\mathsf{T}}e^{C^{\mathsf{T}}_kt}\mathsf{d}  t,\notag
\end{align}
Note that
\begin{align}\label{cg}
	&\mathsf{diag}(U^{-1},I_m)(C_k,\bar{G}_k)\mathsf{diag}(U,I_{2m+p})\notag\\
	&=\left(\left(\begin{array}{ccc}
		A_1 & A_2 & B_1   \\ 
		0 & A_3 & 0   \\ 
		* & * & * \\
	\end{array}\right),\left(\begin{array}{cc}
		D_1 & 0  \\
		0   & 0  \\
		* & I_m   \\ 
	\end{array}\right)\right).
\end{align}
Then for any $\eta\in\text{Im}(A|B_0)\oplus\mathbb{R}^n$ with $\eta\neq0$, it follows that there is some $c>0$ such that
\begin{align}
	&\liminf_{k\to\infty}\gamma_k^{-2}\eta^{\mathsf{T}}\mathsf{diag}(U^{-1},I_m)(N(k+1)-\notag\\
	&\quad\quad\quad\quad N(k))\mathsf{diag}(U^{-1},I_m)^{\mathsf{T}}\eta>c.
\end{align}
Otherwise, there is a sequence $\{k_i,i\in\mathbb{N}\}$ such that 
\begin{align}\label{k9}
	0&=\lim_{{k_i}\to\infty}\int_{0}^{1}(1-t)\eta^{*\mathsf{T}}e^{C_{k_i}t}\bar{G}_{k_i}\bar{G}_{k_i}^{\mathsf{T}} e^{C_{k_i}^{\mathsf{T}} t}\eta^*\mathsf{d}   t\notag\\
	&=\int_{0}^{1}(1-t)\eta^{*\mathsf{T}}e^{Ct}\bar{G}\bar{G}^{\mathsf{T}} e^{C^{\mathsf{T}} t}\eta^*\mathsf{d}   t
\end{align}
where $\eta^*=\mathsf{diag}(U^{-1},I_m)^{\mathsf{T}}\eta$ and $(C,\bar{G})=\lim_{{k_i}\to\infty}(C_{k_i},\bar{G}_{k_i})$. Note that $(C,\bar{G})$ still has the form (\ref{cg}), so $\eta^*\in\text{Im}(C|\bar{G})$, which is contradictory to (\ref{k9}). Hence, (\ref{n}) is true.

Now, we have
\begin{align}
	&\lambda_{\min}(\int_{0}^{k}U_1\varphi(s)\varphi^{\mathsf{T}}(s)U_1^{\mathsf{T}}\mathsf{d}  s)\notag\\
	&=\lambda_{\min}(U_1H(k)U_1^{\mathsf{T}})\notag\\
	&\ge \lambda_{\min}(U_1N(k)U_1^{\mathsf{T}}+U_1M(k)U_1^{\mathsf{T}})\notag\\
	&\ge\sum\limits_{i=0}^{k-1}\lambda_{\min}(U_1(N(i+1)-N(i))U_1^{\mathsf{T}})-\|M(k)\|\|U_1\|^2\notag\\
	&\ge\sum\limits_{i=0}^{k-1}c\gamma_{i}^2-\|M(k)\|\|U_1\|^2.\notag
\end{align}
Therefore, by Lemma \ref{le5}, for some $c_0>0$ and sufficiently large $k\in\mathbb{N}$, we have
\begin{equation}\label{phi2}
	\lambda_{\min}(\int_{0}^{k}U_1^{\mathsf{T}}\varphi(s)\varphi^{\mathsf{T}}(s)U_1^{\mathsf{T}}\mathsf{d} s)\ge c_0k^{3/5}>k^{1/2}.
\end{equation}

Next, it is easy to verify that
\begin{align}
	&\lambda_{\max}\big(U_2(\int_{0}^{k}\varphi(s)\varphi^{\mathsf{T}}(s)\mathsf{d} s)U_2^{\mathsf{T}}\big)\notag\\
	&=\lambda_{\max}\big(\int_{0}^{k}\mathsf{diag}(0_{n_1},x_2^{\mathsf{T}}x_2,0_m)\mathsf{d} s\big)\notag\\
	&=O(1),\notag
\end{align}
which completes the proof.

\vspace{5mm}
\textbf{Proof of Lemma \ref{le7}}:
Let us denote $$\tilde{\theta}^{\mathsf{T}}(k)=[\tilde{\theta}_1^{\mathsf{T}}(k),\tilde{\theta}_2^{\mathsf{T}}(k)] \text{ with } \tilde{\theta}(k)=\hat{\theta}(k)-\theta,$$
\begin{align}
	\left(\begin{array}{cc}
		\Phi_{11}(k) & \Phi_{12}(k) \\ 
		\Phi_{21}(k) & \Phi_{22}(k) \\ 
	\end{array}\right)&\triangleq P^{-1}(k-1)\notag\\
	&=I_{m+n}+\int_{0}^{k-1}a(s)\varphi(s)\varphi^{\mathsf{T}}(s)\mathsf{d}s,\notag
\end{align}
and $e_i\in\mathbb{R}^{n}$ whose $i\text{-th}$  element is $1$ and the other elements are $0$.

Without loss of generality, we assume that $$\lim_{k\to\infty}\tilde{\theta}_1(k) e_i\neq 0,$$
otherwise we can delete the i-th colunm of $\tilde{\theta}(k)$ since it is already convergent to zero. 

It is easy to verify that
\begin{align}
	&\frac{  e_i^{\mathsf{T}}\tilde{\theta}_2^{\mathsf{T}}(k)\Phi_{22}(k)\tilde{\theta}_2(k)  e_i}{ e_i^{\mathsf{T}}\tilde{\theta}_1^{\mathsf{T}}(k)\Phi_{11}(k)\tilde{\theta}_1(k) e_i}\notag\\
	&=O(\frac{\lambda_{\max}(\int_{0}^{k-1}a(s)\varphi_2(s)\varphi_2^{\mathsf{T}}(s)ds)}{\lambda_{\min}\int_{0}^{k-1}(a(s)\varphi_1(s)\varphi_1^{\mathsf{T}}(s)ds)})\notag\\
	&=O(\frac{\lambda_{\max}(\int_{0}^{k-1}\varphi_2(s)\varphi_2^{\mathsf{T}}(s)ds)}{\lambda_{\min}\int_{0}^{k-1}(a(k)\varphi_1(s)\varphi_1^{\mathsf{T}}(s)ds)})\notag\\
	&=o(1).\notag
\end{align}
By Cauchy-Schwarz inequality, it is easy to verify that
\begin{align}
	&e_i^{\mathsf{T}}\tilde{\theta}_1^{\mathsf{T}}(k)\Phi_{12}(k)\tilde{\theta}_2(k)e_i\notag\\
	&= e_i^{\mathsf{T}}\tilde{\theta}_1^{\mathsf{T}}(k)\Phi_{12}(k)\tilde{\theta}_2(k)  e_i\notag\\
	&\le( e_i^{\mathsf{T}}\tilde{\theta}_1^{\mathsf{T}}(k)\Phi_{11}(k)\tilde{\theta}_1(k) e_i)^{\frac{1}{2}}(  e_i^{\mathsf{T}}\tilde{\theta}_2^{\mathsf{T}}(k)\Phi_{22}(k)\tilde{\theta}_2(k)  e_i)^{\frac{1}{2}}\notag\\
	&=( e_i^{\mathsf{T}}\tilde{\theta}_1^{\mathsf{T}}(k)\Phi_{11}(k)\tilde{\theta}_1(k) e_i)\frac{( e_i^{\mathsf{T}}\tilde{\theta}_2^{\mathsf{T}}(k)\Phi_{22}(k)\tilde{\theta}_2(k) e_i)^{\frac{1}{2}}}{( e_i^{\mathsf{T}}\tilde{\theta}_1^{\mathsf{T}}(k)\Phi_{11}(k)\tilde{\theta}_1(k) e_i)^{\frac{1}{2}}}\notag\\
	&=( e_i^{\mathsf{T}}\tilde{\theta}_1^{\mathsf{T}}(k)\Phi_{11}(k)\tilde{\theta}_1(k) e_i)\cdot o(1)\notag\\
	&=o(e_i^{\mathsf{T}}\tilde{\theta}_1^{\mathsf{T}}(k)\Phi_{11}(k)\tilde{\theta}_1(k)e_i).\notag
\end{align}
Then, it follows that
\begin{align}
	&e_i^{\mathsf{T}}\tilde{\theta}^{\mathsf{T}}(k)P^{-1}(k-1)\tilde{\theta}(k)e_i\notag\\
	&=e_i^{\mathsf{T}}\tilde{\theta}_1^{\mathsf{T}}(k)\Phi_{11}(k)\tilde{\theta}_1(k)e_i+e_i^{\mathsf{T}}\tilde{\theta}_1^{\mathsf{T}}(k)\Phi_{12}(k)\tilde{\theta}_2(k)e_i\notag\\
	&\quad+e_i^{\mathsf{T}}\tilde{\theta}_2^{\mathsf{T}}(k)\Phi_{21}(k)\tilde{\theta}_1(k)e_i+e_i^{\mathsf{T}}\tilde{\theta}_2^{\mathsf{T}}(k)\Phi_{22}(k)\tilde{\theta}_2(k)e_i \notag\\
	&=\!e_i^{\mathsf{T}}\tilde{\theta}_1^{\mathsf{T}}(k)\Phi_{11}(k)\tilde{\theta}_1(k)e_i(1\!+\!o(1))\!+\!O(\|\Phi_{22}(k)\|),\notag
\end{align}

From lemma \ref{le3} 4), we have 
$$ \|\tilde{\theta}^{\mathsf{T}}(k)P^{-1}(k-1)\tilde{\theta}(k)\|=O(1),$$
then it follows that
$$e_i^{\mathsf{T}}\tilde{\theta}_1^{\mathsf{T}}(k)\Phi_{11}(k)\tilde{\theta}_1(k)e_i=O(\|\Phi_{22}(k)\|).$$
Therefore, we have
\begin{align}
	\text{tr}(\tilde{\theta}_1^{\mathsf{T}}(k)\tilde{\theta}_1(k))&=\sum\limits_{i=1}^ne_i^{\mathsf{T}}\tilde{\theta}_1^{\mathsf{T}}(k)\tilde{\theta}_1(k)e_i\notag\\
	&=O(\frac{\lambda_{\max}(\Phi_{22}(k))}{\lambda_{\min}(\Phi_{11}(k))})\notag\\
	&=O(\frac{\lambda_{\max}\int_{0}^{k-1}(\varphi_2(s)\varphi_2^{\mathsf{T}}(s)ds}{a(k)\lambda_{\min}\int_{0}^{k-1}(\varphi_1(s)\varphi_1^{\mathsf{T}}(s)ds)})\notag\\
	&=o(1).\notag
\end{align}
Hence, the proof is completed.

\bibliographystyle{ieeetr}

\begin{thebibliography}{10}
	
	\bibitem{a3}R. E. Kalman, ``Design of a self-optimizing control system," {\it  Trans. ASME.}, 80, pp. 468-477, 1958.
	\bibitem{b1}R. Bellman, {\it Adaptive Control Processes - A Guided Tour}, Princeton University Press, New Jersey, 1961.
	\bibitem{a4}K. J. Astrom, ``Adaptive control around 1960," {\it  IEEE Control Systems}, 16, pp. 44-49, 1996.
	\bibitem{b2}H. F. Chen and L. Guo, {\it Identification and Stochastic Adaptive Control}, Boston, MA: Birkhauser, 1991.
	\bibitem{b3}K. J. Astrom and B. Wittenmark, {\it Adaptive Control}, Dover Publications, Mineola, N.Y., 2008.
	\bibitem{a5}K. J. Astrom and B. Wittenmark, ``On self-tuning regulators," {\it Automatica}, 9, pp. 185-199, 1973.
	\bibitem{a6}L. Guo and H. F. Chen, ``The Astrom-Wittenmark self-tuning regulator revisited and ELS-based adaptive trackers," {\it IEEE Trans. Automat. Contr.}, 36, pp. 802-812, 1991.
	\bibitem{a7}L. Guo, ``Convergence and logarithm laws of self-tuning regulators," {\it Automatica}, 31, pp. 435–450, 1995.
	\bibitem{a8}H. F. Chen and L. Guo, ``Optimal adaptive control and consistent parameter estimates for ARMAX model with quadratic cost," {\it SIAM J. Control Optim.}, 25, pp. 845-867, 1987.
	\bibitem{a10}H. F. Chen and L. Guo, ``Convergence rate of least-squares identification and adaptive control for stochastic system," {\it International Journal of Control}, 44, pp. 1459-1476, 1986,.
	
	\bibitem{a1}L. Guo, “Self-convergence of weighted least-squares with applications
	to stochastic adaptive control,” {\it IEEE Trans. Automat. Contr.}, 41, pp. 79–89, 1996.
	
	\bibitem{a2}T. E. Duncan, L. Guo and B. Pasik-Duncan, ``Adaptive continuous-time linear quadratic Gaussian control," {\it IEEE Trans. Automat. Contr.}, vol. 44, no. 9, pp. 1653-1662, 1999.
	
	\bibitem{d1}J. G. Ziegler and N. B. Nichols, “Optimum settings for automatic controllers,” {\it Trans. ASME}, 64, pp. 759–768, 1942.
	\bibitem{d2} S. J. Bradtke, “Reinforcement learning applied to linear quadratic regulation,” {\it Proc. Adv. Neural Inf. Process. Syst.}, pp. 295–302, 1993.
	\bibitem{d3} R. E. Skelton and G. Shi, “The data-based LQG control problem,”  {\it Proc.
		IEEE Conf. Decis. Control}, pp. 1447–1452, 1994.
	\bibitem{d4} U. S. Park and M. Ikeda, “Stability analysis and control design of LTI
	discrete-time systems by the direct use of time series data,” {\it Automatica},
	45, pp. 1265–1271, 2009.
	\bibitem{d5} Z. Wang and D. Liu, “Data-based controllability and observability analysis
	of linear discrete-time systems,” {\it IEEE Trans. Neural Netw.}, 22,
	pp. 2388–2392, 2011.
	\bibitem{d6}D. Liu, P. Yan, and Q. Wei, “Data-based analysis of discrete-time linear
	systems in noisy environment: Controllability and observability,” {\it Inf. Sci.},
	288, pp. 314–329, 2014.
	\bibitem{d7}H. Niu, H. Gao, and Z. Wang, “A data-driven controllability measure
	for linear discrete-time systems,” {\it Proc. IEEE 6th Data Driven Control
		Learn. Syst. Conf.}, pp. 455–466, 2017.
	%\bibitem{d8}B. Zhou, Z. Wang, Y. Zhai, and H. Yuan, “Data-driven analysis methods for controllability and observability of a class of discrete LTI systems with delays,” in Proc. IEEE 7th Data Driven Control Learn. Syst. Conf., May 2018, pp. 380–384.
	%\bibitem{d9}T. M. Maupong, J. C. Mayo-Maldonado, and P. Rapisarda, “On Lyapunov functions and data-driven dissipativity,” IFAC-PapersOnLine, vol. 50, no. 1, pp. 7783–7788, 2017.
	%\bibitem{d10}A. Romer, J. M. Montenbruck, and F. Allgöwer, “Determining dissipation inequalities from input-output samples,” IFAC-PapersOnLine, vol. 50, no. 1, pp. 7789–7794, 2017.
	%\bibitem{d11} J. C. Willems, P. Rapisarda, I. Markovsky, and B. L. M. De Moor, “A note on persistency of excitation,” Syst. Control Lett., vol. 54, no. 4, pp. 325–329, 2005.
	%\bibitem{d12} I. Markovsky and P. Rapisarda, “Data-driven simulation and control,” Int. J. Control, vol. 81, no. 12, pp. 1946–1959, 2008.
	%\bibitem{d13}A. Romer, J. Berberich, J. Köhler, and F. Allgöwer, “One-shot verification of dissipativity properties from input-output data,” IEEE Control Syst. Lett., vol. 3, no. 3, pp. 709–714, Jul. 2019.
	%\bibitem{d14}C. De Persis and P. Tesi, "Formulas for Data-Driven Control: Stabilization, Optimality, and Robustness," in IEEE Transactions on Automatic Control, vol. 65, no. 3, pp. 909-924, March 2020.
	%\bibitem{d15}H. J. van Waarde, J. Eising, H. L. Trentelman and M. K. Camlibel, "Data Informativity: A New Perspective on Data-Driven Analysis and Control," in IEEE Transactions on Automatic Control, vol. 65, no. 11, pp. 4753-4768, Nov. 2020.
	
	
	\bibitem{r1}B. Kiumarsi, F. L. Lewis, and Z. P. Jiang, “$H_{\infty}$ control of linear discretetime systems: Off-policy reinforcement learning,” {\it Automatica}, 78,
	pp. 144–152, 2017.
	\bibitem{r2}H. Modares and F. L. Lewis, “Linear quadratic tracking control of
	partially-unknown continuous-time systems using reinforcement learning,” {\it IEEE Trans. Automat. Control}, 59, pp. 3051–3056, 2014.
	\bibitem{r3}M. Fazel, R. Ge, S. M. Kakade, and M. Mesbahi, “Global convergence of
	policy gradient methods for the linear quadratic regulator,” {\it  Proc. Int.
		Conf. Mach. Learn.}, pp. 1467–1476, 2018.
	\bibitem{r4}T. Bian, Y. Jiang, and Z. P. Jiang, “Adaptive dynamic programming for
	stochastic systems with state and control dependent noise,” {\it IEEE Trans.
		Automat. Control}, 61, pp. 4170–4175, 2016.
	%\bibitem{r5}Y. Abbasi-Yadkori, N. Lazic, and C. Szepesvari, “Model-free linear quadratic control via reduction to expert prediction,” in Proc. 22nd Int. Conf. Artif. Intell. Statist., 2019, pp. 3108–3117.
	%\bibitem{r6} S. Tu and B. Recht, “Least-squares temporal difference learning for the linear quadratic regulator,” in Proc. Int. Conf. Mach. Learn., 2018, pp. 5005–5014.
	%\bibitem{r7}S. Dean, H. Mania, N. Matni, B. Recht, and S. Tu, “On the sample complexity of the linear quadratic regulator,” Found. Comput. Math., vol. 20, no. 4, pp. 633–679, 2019.
	%\bibitem{r8}F. A. Yaghmaie, F. Gustafsson and L. Ljung, "Linear Quadratic Control Using Model-Free Reinforcement Learning," in IEEE Transactions on Automatic Control, vol. 68, no. 2, pp. 737-752, Feb. 2023.
	
	\bibitem{r9}S. J. Bradtke, B. E. Ydestie and A. G. Barto, “Adaptive linear quadratic control using policy iteration,” {\it Proc. Amer. Control Conf.}, pp. 3475–3476, 1994.
	\bibitem{r11} P. Werbos, “Neural networks for control and system identification,” {\it Proc. Conf. Decis. Control}, pp. 260–265, 1989.
	%\bibitem{r12} S. A. A. Rizvi and Z. Lin, “Output feedback Q-learning control for the discrete-time linear quadratic regulator problem,” IEEE Trans. Neural Netw. Learn., vol. 30, no. 5, pp. 1523–1536, May 2019.
	%\bibitem{r13}B. Luo, D. Liu, and H.-N. Wu, “Adaptive constrained optimal control design for data-based nonlinear discrete-time systems with critic-only structure,” IEEE Trans. Neural Netw. Learn., vol. 29, no. 6, pp. 2099–2111, Jun. 2018.
	%\bibitem{r14} B. Kiumarsia, B. AlQaudi, H. Modaresa, F. L. Lewis, and D. S. Levinee, “Optimal control using adaptive resonance theory and Q-learning,” Neurocomputing, vol. 361, pp. 119–125, 2019.
	%\bibitem{r15}L. C. Baird, “Reinforcement learning in continuous time: Advantage updating,” in Proc. IEEE Int. Conf. Neural Netw., 1994, pp. 2448–2453.
	\bibitem{r16}T. Bian and Z. P. Jiang, “Reinforcement learning for linear continuous-time systems: An incremental learning approach,” {\it IEEE/CAA J. Automatica Sinica}, 6, , pp. 433–440, 2019.
	\bibitem{r17}Y. Jiang and Z. P. Jiang, “Global adaptive dynamic programming for
	continuous-time nonlinear systems,” {\it IEEE Trans. Autom. Control}, 60, pp. 2917–2929, 2015.
	\bibitem{r18} J. Lee and R. S. Sutton, “Policy iterations for reinforcement learning
	problems in continuous time and space-fundamental theory and methods,”
	{\it Automatica}, 126, Art no. 109421, 2021.
	\bibitem{r19}D. Vrabie, O. Pastravanu, M. Abu-Khalaf, and F. L. Lewis, “Adaptive optimal control for continuous-time linear systems based on policy iteration,”
	{\it Automatica}, 45, pp. 477–484, 2009.
	\bibitem{r20} H. Modares, F. L. Lewis, and Z.-P. Jiang, “Optimal output-feedback control
	of unknown continuous-time linear systems using off-policy reinforcement learning,” {\it  IEEE Trans. Cybern.}, 46, pp. 2401–2410, 2016.
	%\bibitem{r21}N. Li, X. Li, J. Peng and Z. Q. Xu, "Stochastic Linear Quadratic Optimal Control Problem: A Reinforcement Learning Method," in IEEE Transactions on Automatic Control, vol. 67, no. 9, pp. 5009-5016, Sept. 2022.
	\bibitem{Liu2024} N. Liu and L. Guo, ``Adaptive stabilization of noncooperative
	stochastic differential games," {\it SIAM SIAM J. Control Optim.}, 62, pp. 1317–1342, 2024.
	
	%\bibitem{g1}N. Liu and L. Guo, "Stochastic Adaptive Linear Quadratic Differential Games," in IEEE Transactions on %Automatic Control, vol. 69, no. 2, pp. 1066-1073, Feb. 2024,
	\bibitem{b7} Y. Chow and H. Teicher, {\it Probability Theory: Independence, Interchangeability, Martingales}, Springer-Verlag, 2008.
	\bibitem{b6} P. E. Caines, {\it Linear Stochastic Systems}, Wiley, New York, 1988.
	
\end{thebibliography}

\end{document}